\documentclass[12pt]{amsart}
\usepackage[]{amsmath, amsthm, amsfonts, amssymb, epsfig}
\newcommand\dee{\partial}
\newcommand\Om{\Omega}
\newcommand\Obar{\overline{\Omega}}

\numberwithin{equation}{section}

\begin{document}

\title[Dirichlet and Neumann
problems in quadrature domains]
{The Dirichlet and Neumann and Dirichlet-to-Neumann
problems in quadrature, double quadrature, and non-quadrature
domains}
\author[S.~R.~Bell]
{Steven R.~Bell}

\address[]{Mathematics Department, Purdue University, West Lafayette,
IN  47907}
\email{bell@math.purdue.edu}
\thanks{Research supported by the NSF Analysis and Cyber-enabled
Discovery and Innovation programs, grant DMS~1001701}


\begin{abstract}
We demonstrate that solving the classical problems mentioned
in the title on quadrature domains when the given boundary data
is rational is as simple as the method of partial fractions.  A
by-product of our considerations will be a simple proof that
the Dirichlet-to-Neumann map on a double quadrature domain sends
rational functions on the boundary to rational functions
on the boundary. The results extend to more general domains
if rational functions are replaced by the class of functions
on the boundary that extend meromorphically to the double.
\end{abstract}

\maketitle

\theoremstyle{plain}

\newtheorem {thm}{Theorem}[section]
\newtheorem {lem}[thm]{Lemma}

\hyphenation{bi-hol-o-mor-phic}
\hyphenation{hol-o-mor-phic}

\section{Introduction}
\label{sec1}

It has recently come to light that double
quadrature domains in the plane exist in great
abundance and that they can be viewed as replacements
for the unit disc when it comes to questions of computational
complexity in conformal mapping and potential theory. They
are especially useful in the multiply connected setting.
The ``improved Riemann mapping theorem''
described in \cite{BGS} and further expounded
upon in \cite{B8} allows one to map domains in the
plane, even multiply connected domains, to {\it nearby\/}
double quadrature domains, thus providing the means
to pull back objects on double quadrature
domains to the original domain. Double quadrature
domains share many of the beautiful and simple
properties of the unit disc. The purpose of this
paper is to explain methods for solving some of
the classical problems of potential theory in
quadrature domains that are every bit as simple
as similar problems on the unit disc. In fact,
the methods will be shown to be analogous to
the method of partial fractions from freshman
calculus.

In this paper, we will consider problems in
potential theory with rational boundary data in two
special types of domains. The first type will be
bounded quadrature domains with respect to area
measure without cusps in the boundary. We shall
call such domains {\it area quadrature domains}.
The second type will be area quadrature domains
which are also quadrature domains with respect
to boundary arc length. We shall call such domains
{\it double quadrature domains}. We refer the
reader to \cite{BGS} for the precise definitions
of these domains and for a summary of their basic
properties. The theory of area quadrature domains was
pioneered by Aharonov and Shapiro in \cite{AS} in
the simply connected setting and by Gustafsson
\cite{G} in the multiply connected setting. Quadrature
domains with respect to boundary arc length were
studied by Shapiro and Ullemar in the simply
connected setting in \cite{SU} and by Gustafsson
in \cite{G2} in the multiply connected setting.
Good references for information about quadrature
domains, their usefulness, and the history of the
subject are the book \cite{EGKP}, the papers
\cite{GS} and \cite{C} therein, and the book \cite{S}.

We now list some of the key properties of area and
double quadrature domains that we will need in
what follows. To begin, assume that $\Om$ is
an area quadrature domain. Then $\Om$ has a
boundary consisting of finitely many non-intersecting
$C^\infty$ smooth real analytic curves which are, in
fact, real algebraic. Gustafsson \cite{G} (after
Aharonov and Shapiro \cite{AS} in the simply connected
case) showed that the Schwarz function $S(z)$ associated
to $\Om$ extends meromorphically to the double
of $\Om$. (A meromorphic function $h$ on $\Om$ that
extends continuously up to the boundary extends
meromorphically to the double if and only if there
is a meromorphic function $H$ on $\Om$ that also
extends continuously to the boundary such that
$h=\overline{H}$ on the boundary.) Since
\begin{equation}
\label{schwarz}
S(z)=\bar z
\end{equation}
on the boundary, it follows that $z$ extends
meromorphically to the double (by defining the
extension to be $\overline{S(z)}$ on the backside)
and $S(z)$ extends meromorphically to the double
(by defining the extension to be $\bar z$ on
the backside). Gustafsson showed that the meromorphic
extensions of the two functions $z$ and $S(z)$ to
the double form a primitive pair for the double,
meaning that they generate the field of meromorphic
functions on the double. (See Farkas and Kra \cite{FK}
for the basic facts about primitive pairs and the
field of meromorphic functions on the double.)
Identity~(\ref{schwarz}) allows us to see that a
rational function $R(z,\bar z)$ of $z$ and $\bar z$ is
equal to $R(z,S(z))$ on the boundary and this yields an
extension of the rational function to the double as a meromorphic
function. Conversely, if $G$ is a meromorphic function
on the double, then $G$ is a rational combination of
the extensions of $z$ and $S(z)$. Since $S(z)=\bar z$
on the boundary, we see that the restriction of $G$
to the boundary is rational.
We have seen, therefore, that the field of rational
functions $R(z,\bar z)$ on the boundary is precisely
the set of meromorphic functions on the double restricted
to the boundary. Consequently, if the rational
function does not blow up on the boundary, then
it is $C^\infty$ smooth on the boundary. Let
${\mathcal R}(b\Om)$ denote the set of rational functions
on the boundary without singularities on the boundary, and let
${\mathcal R}(\Om)$ denote the set of meromorphic functions
on $\Om$ obtained by extending functions in
${\mathcal R}(b\Om)$ to $\Om$ via the formula
$R(z,\bar z)=R(z,S(z))$.

The Szeg\H o kernel $S(z,w)$ associated to the area
quadrature domain $\Om$ extends holomorphically in
$z$ and anti-holomorphically in $w$ to an open set
containing $\Obar\times\Obar$ minus the boundary
diagonal. The Garabedian kernel $L(z,w)$ 
extends holomorphically in $z$ and holomorphically
in $w$ to an open set containing $\Obar\times\Obar$
minus the diagonal. It has a simple pole in the $z$
variable at $z=w$ when $w\in\Om$ is held fixed. The
residue in $z$ at $w$ is $1/2\pi$. The Szeg\H o kernel
and Garabedian kernel are non-vanishing in simply
connected domains, but on an $n$-connected domain,
the Szeg\H o kernel $S(z,w)$ has $n-1$ zeroes in
$z$ on $\Om$ for each fixed $w$ in $\Om$. The
Garabedian kernel $L(z,w)$, however, is non-zero if
$z\in\Obar$ and $w\in\Om$ with $z\ne w$ even in the
multiply connected case. If
$a\in\Om$, then neither $S(z,a)$ nor $L(z,a)$
vanish for $z$ in the boundary. See \cite{B1} for
proofs of all these facts in the spirit of this paper.

Let $S^0(z,w)$ denote $S(z,w)$ and let $S^m(z,w)$
denote $(\dee/\dee\bar w)^mS(z,w)$. Similarly, 
let $L^0(z,w)$ denote $L(z,w)$ and let $L^m(z,w)$
denote $(\dee/\dee w)^mL(z,w)$.
The {\it Szeg\H o span\/} is the complex linear
span $\mathcal S$ of all functions $h(z)$ of the form
$h(z)=S^m(z,a)$ as $a$ ranges over $\Om$ and $m$
ranges over all non-negative integers.
The {\it Garabedian span\/} is the complex linear
span $\mathcal L$ of all functions $H(z)$ of the form
$H(z)=L^m(z,a)$ as $a$ ranges over $\Om$ and $m$
ranges over all non-negative integers.
The {\it Szeg\H o plus Garabedian span\/} is the set
$\mathcal S +\mathcal L$ of all sums $h+H$ where
$h$ is in the Szeg\H o span and $H$ is in the Garabedian
span. We will often shorten our notation by writing
$S_a^m(z)=S^m(z,a)$ and $L_a^m(z)=L^m(z,a)$, and
we emphasize here that the unadorned $S(z)$ will
always stand for the Schwarz function. Note that,
because $L(z,a)$ has a singular part that is 
a non-zero constant times $(z-a)^{-1}$, the singular
part of $L_a^m$ is a non-zero constant times
$(z-a)^{-(m+1)}$.

To have a proper feeling for the Szeg\H o and Garabedian
spans, we mention here that
the Szeg\H o span is a dense subspace of the $L^2$-Hardy
space and the Garabedian span is a dense subspace of the
orthogonal complement to the $L^2$-Hardy space in $L^2(b\Om)$.
Hence $\mathcal S+\mathcal L$ is dense in $L^2$ of the boundary.
Let $A^\infty(\Om)$ denote the space of holomorphic
functions on $\Om$ in $C^\infty(\Obar)$.
The Szeg\H o span is also a dense subspace of
$A^\infty(\Om)$ 
and the Garabedian span is a dense subspace of the
orthogonal complement to the $L^2$-Hardy space in
the topology of $C^\infty(b\Om)$.
Hence $\mathcal S+\mathcal L$ is dense in $C^\infty$
on the boundary.
(See \cite{B0} and \cite{B1} for proofs of these facts
in the more general smooth domain case. The density of
the space of rational functions on the boundary of an
area quadrature domain in $C^\infty$ of the boundary
is also proved in the last section of \cite{B8}.)

Let $T(z)$ denote the complex unit tangent vector
function defined on the boundary of $\Om$ and pointing
in the direction of the standard orientation of the
boundary.  A very important identity at the heart of much
of this paper is the relationship between the Szeg\H o
kernel and the Garabedian kernel,
\begin{equation}
\label{SL}
\overline{S_a(z)}=\frac{1}{i}L_a(z)T(z),
\end{equation}
which holds for $z\in b\Om$ and $a\in\Om$.
We may differentiate this identity with respect to
$a$ to obtain
\begin{equation}
\label{SL2}
\overline{S_a^m(z)}=\frac{1}{i}L_a^m(z)T(z).
\end{equation}

If $\Om$ is a double quadrature domain, then
all the properties above hold plus the property
(proved by Gustafsson in \cite{G2})
that $T(z)$ extends to the double as a meromorphic
function. Consequently, identities~(\ref{SL})
and~(\ref{SL2}) show that the functions $S_a^m$ and
$L_a^m$ also extend to the double for each $a$ in
$\Om$ and $m\ge0$.

With these preliminaries behind us, we can state our
main results. We call the first result the {\it Basic
Decomposition}.

\begin{thm}
\label{main}
Given a point $a$ in an area quadrature domain $\Om$, 
$$S_a{\mathcal R}(b\Om)=\mathcal S+\mathcal L,$$
where the function spaces on the right are understood
to be restricted to the boundary. On a double
quadrature domain,
$${\mathcal R}(b\Om)=\mathcal S+\mathcal L.$$
\end{thm}

We will prove this theorem in the next section,
where it will be seen that the coefficients that
appear in the decompositions are determined by
the principal parts of two meromorphic functions.
We remark here that, because the functions on the
left hand side of the equalities in Theorem~\ref{main}
extend meromorphically to $\Om$ via the
$R(z,\bar z)=R(z,S(z))$ substitution, and the functions
on the right also extend, we may also state the
following result.

\begin{thm}
\label{main2}
Given a point $a$ in an area quadrature domain $\Om$, 
$$S_a{\mathcal R}(\Om)=\mathcal S+\mathcal L.$$
On a double quadrature domain,
$${\mathcal R}(\Om)=\mathcal S+\mathcal L.$$
\end{thm}

We will show how the decomposition in Theorem~\ref{main}
can be used to solve the Dirichlet problem in \S\ref{sec3}.
We will consider the Dirichlet-to-Neumann map in
\S\ref{sec4}, and finally the Neumann problem in
\S\ref{sec6}.

Note that functions in $\mathcal S$ are holomorphic
on $\Om$ and functions in $\mathcal L$ have poles.
Hence, it follows as a corollary to Theorem~\ref{main2}
that the class of holomorphic functions on a double
quadrature domain which extend meromorphically to the
double and which have no singularities on the boundary
is exactly equal to the Szeg\H o span.

Since $\mathcal S + \mathcal L$ is an orthogonal sum,
the following theorem is an easy consequence of
Theorem~\ref{main}.

\begin{thm}
\label{szegoproj}
Let $a$ be a point in an area quadrature domain $\Om$.
The Szeg\H o projection associated to $\Om$ maps
$S_a{\mathcal R}(b\Om)$ onto the Szeg\H o span. If
$\Om$ is a double quadrature domain, it maps
${\mathcal R}(b\Om)$ onto the Szeg\H o span.
\end{thm}

The results of the next section will therefore show
that the Szeg\H o projection of a rational function
can be computed via rather straightforward algebra on
quadrature domains.

Before we start proving the decompositions, we
mention that related results hold
for the Bergman kernel and span. Let $B(z,w)$ denote
the Bergman kernel associated to a bounded area
quadrature domain and let $\Lambda(z,w)$ denote the
complimentary kernel (or conjugate kernel) to the
Bergman kernel (see \cite[p.~134]{B1} for the
definition and basic properties of $\Lambda(z,w)$).
We may define the Bergman span $\mathcal B$ and the
complimentary kernel span $\mathnormal{\Lambda}$
exactly as we defined the Szeg\H o span and
Garabedian span. Let ${\mathcal R}'(\Om)$ denote
the set of meromorphic functions on $\Om$ that are
derivatives of functions in ${\mathcal R}(\Om)$.

\begin{thm}
\label{bergman}
On a simply connected area quadrature domain $\Om$, 
$${\mathcal R}'(\Om)=\mathcal B+\mathnormal{\Lambda}.$$
On a multiply connected area quadrature domain,
${\mathcal R}'(\Om)$ is comprised of the functions
in $\mathcal B+\mathnormal{\Lambda}$ with vanishing periods.
\end{thm}

We will explain in \S\ref{sec5} why a major part
of the proof of Theorem~\ref{bergman} should be
attributed to Bj\"orn Gustafsson.

We remark, that since all non-zero functions in the
space $\mathnormal{\Lambda}$ have poles in $\Om$, a corollary
of Theorem~\ref{bergman} is that, on an area quadrature
domain, a function in ${\mathcal R}'(\Om)$ without
singularities in $\Om$ must be in the Bergman span.
This result will allow us to characterize the image
of the rational functions under the Dirichlet-to-Neumann
map of an area quadrature domain.

Call the map that takes Dirichlet problem boundary
data to the normal derivative of the solution to the
Dirichlet problem the D-to-N map. 

\begin{thm}
\label{d-to-n}
On an area quadrature domain, the D-to-N map takes
${\mathcal R}(b\Om)$ into
$\mathcal B T+ \overline{\mathcal B T}$.
On a double quadrature domain,
$\mathcal B T+ \overline{\mathcal B T}$ is contained
in ${\mathcal R}(b\Om)$ and so the D-to-N map takes
takes ${\mathcal R}(b\Om)$ into itself.
\end{thm}

More can be said in case the quadrature domains are
simply connected.

\begin{thm}
\label{d-to-n-SC}
On a simply connected area quadrature domain, the D-to-N map takes
${\mathcal R}(b\Om)$ onto
$\mathcal B T+ \overline{\mathcal B T}$.
\end{thm}

We will show that, on an area quadrature domain, the
decomposition $\mathcal B T+ \overline{\mathcal B T}$
in Theorem~\ref{d-to-n} uniquely
determines the functions in the Bergman span appearing
in the sum, as made precise in the following
theorem.

\begin{thm}
\label{uniqueness}
On an area quadrature domain, functions in
$\mathcal B T+ \overline{\mathcal B T}$
are represented as
$\kappa_1 T+\overline{\kappa_2 T}$
by uniquely determined elements $\kappa_1$ and
$\kappa_2$ in the Bergman span.
\end{thm}

Results like this will allow us to consider a one-sided inverse
to the D-to-N map in the setting of Theorem~\ref{d-to-n-SC} that
is defined in rather explicit terms (see
Theorem~\ref{inverse} in \S\ref{sec5}).

\section{The basic decomposition in area quadrature domains}
\label{sec2}

Suppose that $\Om$ is a bounded area
quadrature domain and suppose $\phi(z)=R(z,\bar z)$ is a
rational function of $z$ and $\bar z$ without singularities
on $b\Om$. We will now explain how to produce a finite
decomposition of $\phi$ on the boundary in terms of the
Szeg\H o kernel and the Garabedian kernel that can be
thought of as an analogue of a ``partial fractions
decomposition'' on the boundary.  The decomposition will
allow us to solve the Dirichlet problem with rational
boundary data in finite terms.

Pick a point $a$ in $\Om$.  Notice that
$$S_a(z)\phi(z)=S_a(z)R(z,\bar z)= S_a(z)R(z,S(z))$$
on the boundary, and this defines a meromorphic extension
$G$ of $S_a\phi$ to $\Om$. We now subtract off the unique
linear combination $\lambda$ of the functions of the form
$L^m_{b_k}$ to remove the poles of the meromorphic function
$G$ on $\Om$.  We will show that
$$h(z):=G(z) -\lambda(z)$$
is a function in the Szeg\H o span $\mathcal S$.
Indeed, if we pair $h$ with a function $g$ in the
dense subset of the Hardy space consisting of functions
in $A^\infty(\Om)$, and note that functions of the
form $L^m_{b_k}$ are orthogonal to the Hardy space,
we may use the identity $S(z)=\bar z$ and identity~(\ref{SL})
to see that
$$\langle g, h\rangle = \int_{z\in b\Om}
g(z)\,\overline{S_a(z)R(z,S(z))}\ ds
=-i\int_{z\in b\Om} g(z) L_a(z)\,\overline{R(\,\overline{S(z)}\,,\bar z)}
\ dz,$$
and the residue theorem shows that this last integral yields
a finite linear combination of values of $g$ and its
derivatives at finitely many points.  Hence, $h$ is
equal to the linear combination of the functions
$S^m_{a_k}$ which would have the same effect when paired
with $g$.

We have shown that there are finitely many points
$a_n$ and $b_n$ in $\Om$ and positive integers
$N_S$, $M_S$, $N_L$, and $M_L$ such that
\begin{equation}
\label{decomp1}
S_a(z)R(z,\bar z) =
\sum_{n=1}^{N_S}\sum_{m=0}^{M_S} A_{nm} S_{a_n}^m(z)
+\sum_{n=1}^{N_L}\sum_{m=0}^{M_L} B_{nm} L_{b_n}^m(z)
\end{equation}
on the boundary of $\Om$.

We remark here that, just as in the method of
undetermined coefficients, the coefficients and
points in the decomposition (\ref{decomp1}) are
uniquely determined by the principal parts of
meromorphic functions with finitely many poles.
Indeed the coefficients $B_{nm}$ and the points
$b_n$ were chosen so that the principal parts of
the sum $\lambda$ match the principal parts
of $S_a(z)R(z,S(z))$.  If we multiply the decomposition
by $T(z)$ and use identities (\ref{SL}) and (\ref{SL2}),
and note that $R(z,\bar z)=R(\,\overline{S(z)}\,,\bar z)$
on the boundary, we obtain after conjugation
\begin{equation}
\label{decomp2}
L_a(z)\overline{R(\,\overline{S(z)}\,,\bar z)} =
\sum_{n=1}^{N_S}\sum_{m=0}^{M_S} \overline{A_{nm}}\,L_{a_n}^m(z)
+\sum_{n=1}^{N_L}\sum_{m=0}^{M_L} \overline{B_{nm}}\,S_{b_n}^m(z),
\end{equation}
and so we see that the coefficients $A_{nm}$ and
the points $a_n$ are determined by the principal parts of
$L_a(z)\overline{R(\,\overline{S(z)}\,,\bar z)}$ in $\Om$.

We have shown that $S_a$ times a rational function is in
the Szeg\H o plus Garabedian span when restricted to the
boundary.  To finish the proof of Theorem~\ref{main} for area
quadrature domains, we need to show that a function
in  $\mathcal S + \mathcal L$ divided by $S_a$, when restricted
to the boundary, is rational and without singularities on
the boundary. If we divide a function in 
$\mathcal S + \mathcal L$ by $S_a$, we obtain a sum of
functions of the form $S_{a_n}^m/S_a$ and $L_{a_n}^m/S_a$.
Such functions extend to the double of $\Om$ as meromorphic
functions because identities~(\ref{SL}) and~\ref{SL2} show
that they are equal to the conjugates of
$L_{a_n}^m/L_a$ and $S_{a_n}^m/L_a$, respectively, on
the boundary. Since these functions extend to the double and
do not have singularities on the boundary, they are therefore
rational combinations of $z$ and $S(z)$,
which when restricted to the boundary, are rational functions
of $z$ and $\bar z$ without singularities on the boundary.
This completes the proof of the part of Theorem~\ref{main}
about area quadrature domains.

We remark that the space of functions on the boundary
given by $S_a$ times a rational function is easily seen
to be independent of the point $a$ since quotients of
the form $S_a/S_b$ extend meromorphically to the double,
and are therefore rational on the boundary.

On a double quadrature domain, $S_a$ extends to the
double as a meromorphic function and has no singularities
on the boundary (see \cite{BGS}). Therefore
$S_a\mathcal R = \mathcal R$ and the proof of
Theorem~\ref{main} is complete.

\section{Using the decomposition to solve the Dirichlet problem}
\label{sec3}

We first assume that $\Om$ is a {\it simply connected\/} area quadrature
domain. We will now show how the decomposition (\ref{decomp1})
produces a simple and explicit solution to the Dirichlet problem
with rational boundary data $R(z,\bar z)$. Recall that $S_a$ and
$L_a$ are non-vanishing on $\Obar$ and extend holomorphically
past the boundary in case $\Om$ is simply connected.
Notice that if we divide the decomposition
(\ref{decomp1}) by $S_a$, then we decompose our
boundary data $R(z,\bar z)$ into a finite sum of
functions of the form $S_b^m/S_a$ and $L_b^m/S_a$.
The functions $S_b^m/S_a$ extend holomorphically
to $\Om$ and are smooth up to the boundary.  Identities
(\ref{SL}) and (\ref{SL2}) reveal that $L_b^m/S_a$ is
equal to the conjugate of $S_b^m/L_a$ on the boundary and
therefore these functions extend
antiholomorphically to $\Om$ and are smooth up to
the boundary.  (Note that the simple pole of $L_a$ at $a$
gives rise to a zero of the quotient at $a$.)
Consequently, we can read off the conclusion of the
following theorem.

\begin{thm}
\label{thmA}
The solution to the Dirichlet problem on a simply
connected area quadrature domain with rational boundary data
$R(z,\bar z)$ can be read off from the basic decomposition
of $S_a(z)R(z,\bar z)$ and is of the form $h+\overline{H}$
where both $h$ and $H$ are sums of quotients that extend
holomorphically to $\Om$ and meromorphically to the double.
In fact, $h$ is a quotient of the form $\sigma_1/S_a$
and $H$ is a quotient of the form $\sigma_2/L_a$ where
$\sigma_j$, $j=1,2$, are functions in the Szeg\H o span.
\end{thm}

Note that meromorphic functions on the double are generated
by the meromorphic extensions of $z$ and $S(z)$ to the double,
and that $S(z)$ is algebraic. Thus, we have given an alternate
way of looking at Ebenfelt's theorem \cite{E} about the
algebraicity of the solution to the Dirichlet problem
with rational boundary data on simply connected area
quadrature domains (see also \cite{BEKS} for more about
this fascinating subject).

We may repeat much of the same reasoning in case $\Om$
is an $n$-connected area quadrature domain, taking
into account that $S_a$ has $n-1$ zeroes. We may choose
the point $a$ so that the $n-1$ zeroes of $S_a$ are
distinct and simple (see \cite[p.~105]{B1}). Let
$a_1,\dots,a_{n-1}$ denote these zeroes. Let $G(z,w)$
denote the classical Green's function associated to
$\Om$ and write $G_w^0(z)=G(z,w)$. For $m\ge 1$, define
$$G_w^{\,m}(z):=\frac{\dee^m}{\dee w^m}G(z,w)\quad\text{and}\quad
G_w^{\,\bar m}(z):=\frac{\dee^m}{\dee \bar w^m}G(z,w).
$$
Note that, as a function of $z$, $G_b^m(z)$ has a
singular part at $b\in\Om$ that is a non-zero constant
times $1/(z-b)^m$, is harmonic in $z$ on
$\Om-\{b\}$, extends continuously to the boundary
and vanishes on the boundary. Similarly,
$G_b^{\bar m}(z)$ has a
singular part that is a non-zero constant
times the conjugate of $1/(z-b)^m$, is harmonic in $z$ on
$\Om-\{b\}$, extends continuously to the boundary
and vanishes on the boundary.

To solve the Dirichlet problem on $\Om$, given
rational boundary data $R(z,\bar z)$,
as in the simply connected case, the decomposition
(\ref{decomp1}) yields
$$R(z,\bar z)=h+\overline{H}$$
on the boundary, where
$$h=\frac{\sigma_1(z)}{S_a(z)}\qquad\text{ and }\qquad
H=\frac{\sigma_2(z)}{L_a(z)},$$
and where $\sigma_1$ and $\sigma_2$ are in the Szeg\H o
span. Note that $h$ and $H$ extend smoothly to the
boundary, that $H$ is holomorphic on $\Om$ because
$L_a$ is non-vanishing on $\Obar-\{a\}$ and the pole
of $L_a$ creates a zero of $H$ at $a$, but $h$ is
perhaps only meromorphic since it may have simple poles
at some or all of the zeroes of $S_a$. However, by
subtracting off appropriate constants times $G_{a_k}^1$
for each of the zeroes $a_k$ to remove the simple poles,
noting that these functions vanish on the boundary,
we obtain the harmonic extension of our boundary data
to $\Om$ in the form
\begin{equation}
\label{decomp3}
h+\overline{H}+\sum_{k=1}^{n-1}c_kG_{a_k}^1,
\end{equation}
where $h$ and $H$ extend meromorphically to the double,
and hence are rational combinations of $z$ and the
Schwarz function $S(z)$. Consequently, we have an
alternate way to that given in \cite{B8} to see
that the solution to the Dirichlet problem with rational
boundary data is algebraic, modulo an $n-1$ dimensional
subspace.

\section{The Dirichlet-to-Neumann map}
\label{sec4}

We now continue the line of thought of the last section
and consider the implications for the D-to-N map for
rational boundary data on an area quadrature domain
$\Om$. It is shown in \cite[p.~134-135]{B1} that, if $w\in\Om$
is held fixed, the normal derivative $(\dee/\dee n_z)$ of
$G_w^0(z)$ with respect to $z$ is given by
$$\frac{\dee}{\dee n_z}G_w^0(z)=
-iT(z)\frac{\dee}{\dee z}G_w^0(z)
+i\,\overline{T(z)}\,\overline{\frac{\dee}{\dee z}G_w^0(z)}.$$
We remark here that it is also shown there that this
expression can be further manipulated to yield two more
expressions for the same normal derivative,
$$\frac{\dee}{\dee n_z}G_w^0(z)=
-2iT(z)\frac{\dee}{\dee z}G_w^0(z)=
2i\,\overline{T(z)}\,\overline{\frac{\dee}{\dee z}G_w^0(z)}.$$
Although these expressions are shorter and simpler, we will
have reason to prefer the longer form.
Hence, for $m\ge 1$, the normal derivative of $G_w^m(z)$ is
given by
\begin{equation}
\label{normalG}
\frac{\dee}{\dee n_z}G_w^m(z)=
-iT(z)\frac{\dee}{\dee z}G_w^m(z)
+i\,\overline{T(z)}\overline{\frac{\dee}{\dee z}G_w^{\bar m}(z)},
\end{equation}
and the alternative expressions above yield
$$\frac{\dee}{\dee n_z}G_w^m(z)=
-2iT(z)\frac{\dee}{\dee z}G_w^m(z)=
2i\,\overline{T(z)}\,\overline{\frac{\dee}{\dee z}G_w^{\bar m}(z)}.$$

It is shown in \cite[p.~77,~134-135]{B1} that the normal
derivative of the solution to the Dirichlet problem given
by equation (\ref{decomp3}) is
$$-ih'(z)T(z)+i\,\overline{H'(z)T(z)}-i T(z)
\sum_{k=1}^{n-1}c_k
\frac{\dee}{\dee z}
G_{a_k}^1(z)+
i \overline{T(z)}
\sum_{k=1}^{n-1}c_k
\overline{
\frac{\dee}{\dee z}
G_{a_k}^{\bar 1}(z)}.$$
But, for $w\in\Om$, the Bergman kernel $B(z,w)$ is related to
the Green's function via
$$B(z,w)=
-\frac{2}{\pi}\frac{\dee^2}{\dee z\dee\bar w}G(z,w)$$
and the complimentary kernel $\Lambda(z,w)$ to the Bergman kernel
is given by definition as
$$\Lambda(z,w)=
-\frac{2}{\pi}\frac{\dee^2}{\dee z\dee w}G(z,w).$$
(See \cite[p.~134]{B1} for these identities and the basic
properties of $\Lambda(z,w)$.) Consequently,
the last formula for the normal derivative can be rewritten as
\begin{equation}
\label{normal}
-ih'(z)T(z)+i\,\overline{H'(z)T(z)}+
\frac{\pi i}{2}  T(z) \sum_{k=1}^{n-1}c_k \Lambda(z,a_k)
-\frac{\pi i}{2}\,\overline{T(z)}\,\sum_{k=1}^{n-1}c_k
\overline{B(z,a_k)}.
\end{equation}
It is shown in \cite{B5} that, on an area quadrature domain,
if a meromorphic function $g$ extends meromorphically to the
double, then $g'$ also extends meromorphically to the double. It
is shown in \cite{B5} that the Bergman kernel extends
meromorphically to the double on an area quadrature domain.
We will prove momentarily that $\Lambda(z,a_k)$ also extends
meromorphically to the double as a function of $z$. Hence, we
have expressed
the normal derivative of the solution to the Dirichlet problem
as $gT+\overline{GT}$ where $g$ and $G$ are meromorphic functions
on $\Om$ that extend meromorphically to the double, and are
consequently rational combinations of $z$ and $S(z)$. When we
restrict to the boundary, we conclude that $g$ and $G$ are
rational functions of $z$ and $\bar z$ on the boundary. It is
shown in \cite{B5} that, on an area quadrature domain, the
function $T^2$ extends
to the double as a meromorphic function. Hence, the Neumann
boundary data of the solution to the Dirichlet problem with
rational boundary data is a sum of a rational function times
the square root of a rational function plus the conjugate of
such expressions.  On a double quadrature domain, the function
$T$ itself extends meromorphically to the double, and in this
case, we may state that the D-to-N map sends rational functions
of $z$ and $\bar z$ to rational functions of $z$ and $\bar z$.
In the next section, we consider which rational functions appear
in the range of the D-to-N map in this manner.

We now give the promised proof that $\Lambda(z,a)$ extends
meromorphically in $z$ to the double for fixed $a\in\Om$ when
$\Om$ is an area quadrature domain. The Bergman kernel
$B(z,a)$ is related to $\Lambda(z,a)$ via the identity
\begin{equation}
\label{BL}
B(z,a)T(z)=-\overline{\Lambda(z,a)T(z)}
\end{equation}
for $z$ in the boundary.
Since $\Lambda(z,a)$ is
equal to minus the conjugate of the quantity $B(z,a)$ times
$T(z)^2$ on the boundary, and since both of these functions
extend meromorphically to the double, and are therefore
rational functions of $z$ and $S(z)=\bar z$ on the boundary,
it follows that $\Lambda(z,a)$ extends meromorphically to
the double and is equal to a rational combination of $z$
and $\bar z$ on the boundary.

The identity (\ref{BL}) or the simpler expressions for
the normal derivatives of the Green's functions could be
used to simplify (\ref{normal}) to read
\begin{equation}
\label{normala}
-ih'(z)T(z)+i\,\overline{H'(z)T(z)}+
\pi i  T(z) \sum_{k=1}^{n-1}c_k \Lambda(z,a_k),
\end{equation}
or even
\begin{equation}
\label{normalb}
-ih'(z)T(z)+i\,\overline{H'(z)T(z)}
-\pi i\,\overline{T(z)} \sum_{k=1}^{n-1}c_k
\overline{B(z,a_k)},
\end{equation}
but we prefer (\ref{normal}) because the poles
of the $\Lambda(z,a_k)$ terms exactly cancel
the poles of $h'$ the way we have chosen the
coefficients, and therefore the normal derivative
is in fact expressed as $gT+\overline{GT}$ where
$g$ and $G$ are {\it holomorphic\/} functions on
$\Om$ that extend smoothly to the boundary, and
that extend meromorphically to the double. We will
have more to say on this subject later when we
prove Theorem~\ref{d-to-n}.

An interesting consequence of equations~(\ref{normala})
and ~(\ref{normalb}) is that they imply that certain
period matrices are non-singular.

\begin{thm}
\label{periodmatrix}
Suppose that the zeroes $a_1,\dots a_{n-1}$ of the
Szeg\H o kernel associated to a point $a$ in an $n$-connected
area quadrature domain $\Om$ are simple. Then the matrix of
periods associated to the functions $K(z,a_k)$, $k=1,\dots,n-1$
is non-singular. So is the matrix of periods associated to the
functions $\Lambda(z,a_k)$, $k=1,\dots,n-1$.
\end{thm}

We remark that it was proved in \cite{B7} that, if $\Om$
is $n$ connected, then there {\it exist\/} $n-1$ points
$b_1,b_2,\dots,b_{n-1}$ in $\Om$ such that the period
matrix of $K(z,b_k)$ is non-singular. It is also interesting
to note that, because of the way the zeroes of the Szeg\H o
kernel and the periods of the functions in
Theorem~\ref{periodmatrix} transform under conformal
changes of variables, and because smoothly bounded
$n$-connected domains are conformally equivalent to an
area quadrature domain via Gustafsson's theorem \cite{G},
Theorem~\ref{periodmatrix} can be seen to hold for
general bounded smooth $n$-connected domains as well.

To prove Theorem~\ref{periodmatrix}, note that because the
rational functions are dense in $C^\infty(b\Om)$, we may
approximate a harmonic measure function $\omega_k$ which is
harmonic on $\Om$ and equal to one on the $k$-th boundary
curve of the $n-1$ inner boundary curves and equal
to zero on the other boundary curves. The normal derivative of
the solution to the Dirichlet problem with this rational
boundary data is given by equation~(\ref{normala}), and it
can be made as close in $C^\infty(b\Om)$ to the normal
derivative $-iF_k'T$ of $\omega_k$ as desired (see
\cite[p.~87]{B1} for the calculation of this normal
derivative). Since the periods of $F_k'$, $k=1,\dots,n-1$
are well known to be linearly independent, and since
the periods of the functions
given by equation~(\ref{normala}) are linear combinations of
the periods of $\Lambda(z,a_k)$, it follows that the periods
of $\Lambda(z,a_k)$ are independent.  Identity~(\ref{BL})
shows that the periods of $B(z,a_k)$ are just minus the
conjugates of the periods of $\Lambda(z,a_k)$, and so it
follows that the periods of $B(z,a_k)$ are also independent.
This completes the proof of Theorem~\ref{periodmatrix}.

\section{Proof of Theorem~\ref{bergman}}
\label{sec5}
The proof of Theorem~\ref{bergman} follows a similar
pattern to the arguments in the last section and is
motivated by a result of Gustafsson stated as Lemma~{4}
in \cite{B4a}. (In fact, there are two proofs of a closely
related theorem given in \cite{B4a} and a statement
without proof of a converse that is relevant here.
The proof we set out below is a third way of looking
at this problem, and is shorter and simpler than the
arguments given in \cite{B4a}, but it must be said that
the meat of the argument is Gustafsson's idea.)

Suppose that $\Om$ is an area quadrature
domain and that $h$ is a function in ${\mathcal R}(\Om)$.
It follows that $h'$ has only finitely many poles in $\Om$
which are residue free poles of order two or more.
Note that the singular part of $\Lambda(z,a)$ is equal to
a non-zero constant times $(z-a)^{-2}$, and the singular
part of $\Lambda^m(z,a)$ is equal to a non-zero constant
times $(z-a)^{-(m+2)}$. Hence, there is an element
$\lambda$ in the span $\mathnormal{\Lambda}$ that has the
same principal parts at each of the poles of $h'$. But
such a function $\lambda$ is given as the derivative
$(\dee/\dee z)$ of a finite sum of the form
$$\phi(z)=\sum_{n,m}c_{nm}G_{a_n}^m(z).$$
Note that $\phi$ vanishes on the boundary and that $h-\phi$
also has removable singularities at the poles of $h$.
Also note, that because $h$ extends meromorphically to
the double, there is a meromorphic function $H$ on $\Om$ that
extends smoothly to the boundary (and without singularities
on the boundary) such that $h(z)=\overline{H(z)}$ for
$z\in b\Om$.
Now, given $g$ in $A^\infty(\Om)$, we may compute the
inner product $\langle h'-\lambda, g\rangle_{\Om} =$
$$\frac{i}{2}\int_{\Om} (h'-\lambda)\, \bar g\  dz\wedge d\bar z
=\frac{i}{2}\int_{b\Om} (h-\phi) \bar g\ d\bar z 
=\frac{i}{2}\int_{b\Om} h\, \bar g\ d\bar z
=\frac{i}{2}\int_{b\Om} \bar H \, \bar g\ d\bar z,
$$
and the residue theorem implies that this last integral
is equal to the complex conjugate of a finite linear
combination of values of $g$ and its derivatives at
finitely many points in $\Om$. There is an element
$\kappa$ in the Bergman span that has the same effect when
paired with $g$ in the $L^2(\Om)$ inner product. Hence,
since $A^\infty(\Om)$ is dense in the Bergman space,
$h'-\lambda=\kappa$ and we have proved that
$\mathcal{R}'(\Om)\subset \mathcal{B}+\mathnormal{\Lambda}$.

To prove the reverse inclusion, we will need to define
some terminology. Note that we may differentiate identity
(\ref{BL}) with respect to $\bar a$ to obtain the identity
\begin{equation}
\label{BL2}
B^m(z,a)T(z)=-\overline{\Lambda^m(z,a)T(z)}
\end{equation}
for $z$ in the boundary (where the superscript $m$ indicate
derivatives of order $m$ with respect to $\bar a$ in the Bergman
kernel and with respect to $a$ in the $\Lambda$-kernel).
Given a function $g=\kappa+\lambda$ in
$\mathcal{B}+\mathnormal{\Lambda}$, we define the
{\it complimentary function\/} $G$ to $g$ to be the
function gotten by conjugating all the constants in
the linear combination, by changing terms 
of the form $K_b^m$ in $g$ to $\Lambda_b^m$
in the complimentary function, and terms of the form
$\Lambda_b^m$ to $K_b^m$. In this way, we obtain
a function $G$ in $\mathcal{B}+\mathnormal{\Lambda}$
that satisfies
$$g(z)T(z)=-\overline{G(z)T(z)}$$
on the boundary. If $\gamma$ is a curve in the boundary of $\Om$,
then identity (\ref{BL2}) shows that
\begin{equation*}
\int_{\gamma} B^m(z,a)\ dz=
-\int_{\gamma} \overline{\Lambda^m(z,a)}\ d\bar z,
\end{equation*}
and similarly for integrals of $\Lambda^m(z,a)$ by
taking conjugates. Hence, $g$ and $G$ satisfy
\begin{equation}
\label{BL3}
\int_{\gamma} g(z)\ dz=
-\int_{\gamma} \overline{G(z)}\ d\bar z.
\end{equation}
Hence, if a period of $g$ vanishes, then so does the
same period of $G$.

First, we will prove the reverse inclusion in case $\Om$
is a {\it simply connected\/} area quadrature domain. Note
that, in this case, elements
of $\mathcal{B}+\mathnormal{\Lambda}$ have single valued
antiderivatives on $\Om$ which are meromorphic on $\Om$
because all the poles of elements of $\mathnormal{\Lambda}$
are residue free and are of order two or more. Let $h$ be
such an antiderivative,
and write $h' =\kappa +\lambda$ as we did above. Let $G$
be the complimentary function to $h'$ and let $H$ be an
antiderivative of $G$. Let
$\gamma_z$ denote a curve that starts at a boundary point
$b$ and moves along the boundary to another point $z$ in
the boundary. The formula~(\ref{BL3}) holds for the
curve $\gamma_z$ for our complimentary functions $g=h'$
and $G=H'$. Hence, 
$h(z)-h(b)$ is equal to the conjugate of $-(H(z)-H(b))$
on the boundary.  This shows that $h$ extends to the double,
and the proof of the reverse inclusion is complete in the
simply connected case.

We now turn to the multiply connected case.  If
$\kappa+\lambda$ is an element of
$\mathcal{B}+\mathnormal{\Lambda}$
with vanishing periods, then the periods of the complimentary
function are also zero and we obtain two meromorphic
functions $h$ and $H$ as we did above such that
$h'=\kappa+\lambda$ and $h'T=-\overline{H'T}$ on the boundary.
Our task now is to show that we may adjust $h$ and $H$ by
constants in order to make $h=-\overline{H}$ on the
boundary so that we may conclude that $h$ extends to the
double. Choose a point $b$ in the outer boundary of $\Om$
and adjust $h$ and $H$ by subtracting off constants
$h(b)$ and $H(b)$ from $h$ and $H$ so that $h(b)=0$ and
$H(b)=0$. The calculation of the last paragraph shows
$h=-\overline{H}$ on the outer boundary. The key to
seeing that this identity persists on the inner boundaries
is that integrals from $b$ to a point $b'$ on an
inner boundary also agree because of the relationships
between the kernels and the Green's function. Indeed,
if $\gamma$ is a curve in $\Om$ that starts at $b$
and goes to $b'$ and $w\in\Om$ is not in $\gamma$, then
it is well known that
\begin{equation}
\label{green}
\int_\gamma \frac{\dee}{\dee z}G(z,w)\ dz=
-\int_\gamma \frac{\dee}{\dee\bar z}G(z,w)\ d\bar z.
\end{equation}
To see this, note that if $\phi$ is real valued, then
$$\int_\gamma \frac{\dee\phi}{\dee z}\ dz =
\frac12\int_\gamma\left(\frac{\dee\phi}{\dee x}-i
\frac{\dee\phi}{\dee y}\right)(dx+i\,dy),$$
and the real part of this integral is $\phi(b')-\phi(b)$.
Hence, since the Green's function is real and vanishes on
the boundary, the real parts of the two integrals in
(\ref{green}) vanish, and therefore, since they are
conjugates of each other, the identity follows. Now
differentiating (\ref{green}) with respect to $\bar w$
and letting $w=a$ yields that
$$\int_\gamma K^m(z,a)\ dz=
-\int_\gamma \overline{\Lambda^m(z,a)}\ d\bar z.$$
(This is just a reformulation of a well known fact
going back to Bergman and Schiffer about the vanishing
of the $\beta$-periods of the meromorphic differentials
obtained by extending $K_a dz$ to the backside of the
double as the conjugate of $-\Lambda_a dz$.)
To continue, if we now choose such a curve $\gamma$ from $b$ to $b'$
that avoids points where $h'$ and $H'$ have singularities,
this identity shows that
$$h(b')=\int_{\gamma}h'\,dz=
-\int_{\gamma}\overline{H'}\,d\bar z=
-\overline{H(b')}.$$
Consequently, the identity $h=-\overline{H}$ extends
to the inner boundary containing $b'$. We may conclude
that $h$ extends to the double as a meromorphic function.

Theorem~\ref{periodmatrix} yields that if $\Om$
is $n$ connected, then there exist $n-1$ points
$a_1,a_2,\dots,a_{n-1}$ in $\Om$ such that the period
matrix of $K(z,a_k)$ is non-singular. Hence, given an
element $\kappa+\lambda$ of
$\mathcal{B}+\mathnormal{\Lambda}$,
it is possible to subtract off a linear combination of
the functions $K(z,a_k)$ so to make the periods vanish.
Hence, Theorem~\ref{bergman} yields that
$\kappa+\lambda$ is equal to a function in $\mathcal{R}'(\Om)$
modulo a linear combination of $K(z,a_k)$. If we let
$\mathcal{B}_{n-1}$ denote the complex linear span of
the $K(z,a_k)$, then we can state that
$$\mathcal{R}'(\Om)\subset \mathcal{B}+\mathnormal{\Lambda}\subset
\mathcal{R}'(\Om)+\mathcal{B}_{n-1}.$$

Next, we may use Theorem~\ref{bergman} to determine
the image of $\mathcal{R}(b\Om)$ under the Dirichlet-to-Neumann
map in an area quadrature domain $\Om$. Formula~(\ref{normal})
combined with Theorem~\ref{bergman} shows that functions in
the image are of the form
$\kappa_1 T + \lambda_1 T +\overline{\kappa_2 T}
+\overline{\lambda_2 T}$
where $\kappa_1$ and $\kappa_2$ are in $\mathcal B$ and
$\lambda_1$ and $\lambda_2$ are in $\mathnormal{\Lambda}$.
Identity~(\ref{BL2})
shows how to convert the term $\lambda_1 T$ into a term
of the form $\overline{\kappa_3 T}$ where $\kappa_3\in\mathcal B$,
and the term $\overline{\lambda_2T}$ into a term of the
form $\kappa_4 T$ where $\kappa_4\in\mathcal B$.
Hence, when everything is combined, we obtain an expression in
$\mathcal B T+ \overline{\mathcal B T}$.

Next, we show that every element in
$\mathcal B T+ \overline{\mathcal B T}$ is equal to the
normal derivative of a harmonic function with rational
boundary values in case $\Om$ is a simply connected area
quadrature domain. Indeed, given a function
$\psi=\kappa_1 T + \overline{\kappa_2 T}$ in this space,
we may find functions holomorphic functions $h_1$ and $h_2$
on $\Om$ such that $-ih_1'=\kappa_1$ and $-ih_2'=\kappa_2$
where, by Theorem~\ref{bergman}, $h_1$ and $h_2$ are in
$\mathcal{R}(\Om)$. We may now write
$$\psi=-ih_1'T + i\,\overline{h_2'T}.$$
Such a function is the normal derivative of the harmonic
function with rational boundary data $h_1+\overline{h_2}$.

Finally, we need to show that a representation of the
form $\psi=\kappa_1 T + \overline{\kappa_2 T}$ is unique.
If such an expression were equal to zero on
the boundary, then $\kappa_1 T=- \overline{\kappa_2 T}$
and the left hand side of this expression is orthogonal
to holomorphic functions in $L^2(b\Om)$ and the right hand
side is orthogonal to the conjugates of holomorphic
functions in $L^2(b\Om)$. Such functions must be given by
sums of $F_k'T$ where $F_k'$ are the well known holomorphic
functions that arise via $F_k'=2(\dee/\dee z)\omega_k$
where $\omega_k$ are the harmonic measure functions
associated to the $n-1$ inner boundary curves (see
\cite[p.~80]{B1} for a proof of this result). Hence
$\kappa_1=\sum_{k=1}^{n-1}c_kF_k'$. We have just shown
that $\kappa_1 T$ is equal to the normal derivative
of a harmonic function $\phi$ with rational boundary
values. However, $\kappa_1T$ is now also seen to be
the normal derivative of a linear combination $\omega$
of $\omega_k$, $k=1,\dots,n-1$. Hence, $\phi$ and
$\omega$ differ by a constant, and it follows that
the boundary values of $\omega$ are in $\mathcal{R}(b\Om)$.
But the functions in $\mathcal{R}(b\Om)$ are the
boundary values of meromorphic functions of the form
$R(z,S(z))$ and, since $\omega$ vanishes on the outer
boundary, it follows that $\omega\equiv0$, and hence,
that $\phi$ is constant and hence, that $\kappa_1\equiv0$.
Finally, it is an easy exercise to see that if
$$0\equiv \sum_{n=1}^N\sum_{m=0}^M c_{nm}K^m(z,b_m),$$
then all the coefficients $c_{nm}$ must be zero. (Indeed,
such a function would be orthogonal to the Bergman
space, and hence orthogonal to all polynomials. However,
pairing a polynomial with the function in the Bergman
span would yield a non-trivial sum of values and
derivatives of the polynomial at finitely many points
in $\Om$, and it easy to construct a polynomial that
would make this sum non-zero, yielding a contradiction.)
We have shown that the representation of a function
in $\mathcal B T+ \overline{\mathcal B T}$ is uniquely
determined.

The techniques of this section allow us to construct
a one-sided inverse to the D-to-N map on rational
functions in a simply connected area quadrature domain.
Indeed, given a basic term like $K_a^m$ in the Bergman span
we may express a meromorphic
antiderivative of $-ih$ of $K_a^m$ on $\Om$ via
a path integral formula. The proof
of Theorem~\ref{bergman} reveals that $h$ is in
${\mathcal R}(\Om)$. Now, the normal derivative of the
solution to the Dirichlet problem with boundary data
$h\in{\mathcal R}(b\Om)$ is equal to $K_a^mT$. By this
means, we may define a linear transformation $L$ that maps
$K_a^mT$ to the boundary values of $h$.  The same
procedure works for conjugates of terms of the form
$K_a^mT$. Since the representation of functions in
${\mathcal B}T+\overline{{\mathcal B}T}$ is unique,
we obtain the operator $L$ of the following theorem.

\begin{thm}
\label{inverse}
Suppose $\Om$ is a simply connected area quadrature domain.
There is a natural linear transformation $L$ which
maps ${\mathcal B}T+\overline{{\mathcal B}T}$
onto ${\mathcal R}(b\Om)$ such that the D-to-N map
composed with $L$ is the identity.
\end{thm}

\section{The Neumann problem}
\label{sec6}

The Szeg\H o projection can be used to solve the classical
Neumann problem for the Laplacian in planar domains in much
the same way that it was used above to solve the Dirichlet
problem. This process is described in \cite[p.~87]{B1}.
On a double quadrature domain, both the Szeg\H o kernel
$S_a$ and the Garabedian kernel $L_a$ extend to the double
and are therefore rational on the boundary. Also, the
functions $F_j'$ extend to the double (on area quadrature
domains).  If we combine these results with Theorem~20.1
in \cite{B1} and use the fact that the Szeg\H o projection
maps rational functions on the boundary to rational
functions on the boundary, we obtain the following result.

\begin{thm}
\label{neumann}
If $\psi$ is a rational function on the boundary of
a double quadrature domain $\Om$ such that
$\int_{b\Om}\psi\,ds=0$, then the solution to the
Neumann problem with boundary data $\psi$ is equal
to
$$h+\overline{H}+\sum_{k=0}^{n-1} c_k\omega_k$$
where $h$ and $H$ are holomorphic functions on $\Om$
such that $h'$ and $H'$ extend meromorphically to the
double (and are therefore rational on the boundary)
and the $c_k$ are constants.
\end{thm}

\section{Non-quadrature domains}
\label{sec7}

All of the results of this paper carry over to non-quadrature
domains if we define our basic objects differently. Suppose
$\Om$ is a bounded $n$-connected domain with $C^\infty$ smooth
boundary. In this context, let ${\mathcal R}(b\Om)$ denote
the space of $C^\infty$ functions on the boundary that
extend meromorphically to the double of $\Om$, let
${\mathcal R}(\Om)$ denote the space of meromorphic functions
on $\Om$ that have boundary values in
${\mathcal R}(b\Om)$, and let
${\mathcal R}'(\Om)$ denote the space of functions that are
derivatives of functions in
${\mathcal R}(\Om)$. It is proved in \cite{B1b} that there
are two Ahlfors maps $f_1$ and $f_2$ associated to two (rather
generic) points in $\Om$ such that the meromorphic extensions
to the double of $\Om$ form a primitive pair for the double.
Hence, the function spaces just described can all be
expressed in terms of rational functions of $f_1$ and $f_2$.
(Since $\overline{f_j}=1/f_j$ on the boundary $j=1,2$, these
functions conveniently replace the Schwarz function in
many situations.)

The main theorems of the paper in this context can be
stated as follows.

\begin{thm}
\label{maing}
Given a point $a$ in a bounded smooth finitely connected
domain $\Om$,
$$S_a{\mathcal R}(b\Om)=\mathcal S+\mathcal L,$$
where the function spaces on the right are understood
to be restricted to the boundary. Furthermore,
$$S_a{\mathcal R}(\Om)=\mathcal S+\mathcal L.$$
The Szeg\H o projection maps $S_a{\mathcal R}(b\Om)$
onto the Szeg\H o span.
\end{thm}

The Dirichlet problem can be solved for boundary data
in ${\mathcal R}(b\Om)$ by exactly the same methods
we used in \S\ref{sec3}.

The theorem about the Bergman span also generalizes in
a straightforward manner.

\begin{thm}
\label{bergmang}
Suppose that $\Om$ is a bounded smooth finitely connected
domain. If $\Om$ is simply connected, then
$${\mathcal R}'(\Om)=\mathcal B+\mathnormal{\Lambda}.$$
On a multiply connected area quadrature domain,
${\mathcal R}'(\Om)$ is comprised of the functions
in $\mathcal B+\mathnormal{\Lambda}$ with vanishing periods.
In both cases, the D-to-N map takes
${\mathcal R}(b\Om)$ into
$\mathcal B T+ \overline{\mathcal B T}$.
In case $\Om$ is simply connected, this mapping is onto.
Representations of functions $\psi$ in
$\mathcal B T+ \overline{\mathcal B T}$
uniquely determine elements $\kappa_1$ and $\kappa_2$
in the Bergman span so that
$\psi=\kappa_1 T +\overline{\kappa_2 T}$.
\end{thm}

Finally, we remind the reader that we explained in
\S\ref{sec4} why Theorem~\ref{periodmatrix} holds in
general bounded smooth domains. We remark here that
the general result could also be proved from scratch
using the definitions in this section and by repeating
the proof given in \S\ref{sec4} using these definitions.

\end{document}